\documentclass[12pt,reqno]{amsart}
\usepackage{amsmath,amsthm,amscd,amsfonts,amssymb}

\swapnumbers
\numberwithin{equation}{subsection}

\begin{document}
\title{Errata to Local contribution to the Lefschetz fixed
point formula}
\author{Mark Goresky}\thanks{School of Math., Inst. for Adv. Study,
Princeton, N.J.} 
\author{Robert MacPherson}\thanks{School of Math., Inst. for Adv.
Study, Princeton, N.J.}
\begin{abstract}This note contains a correction to the paper, ``Local
contribution to the Lefschetz fixed point formula'', Inv. Math.
{\bf 111} (1993), 1-33. \end{abstract}
\maketitle

This note concerns the paper, ``Local contribution to the Lefschetz
fixed point formula'', Inv. Math. {\bf 111} (1993), 1-33.  We are
grateful to S. Morel for pointing out two errors in this paper.  The
first error occurs in the proof of Lemma 5.10 (which appears in \S 5.11).  
The second error occurs in the proof of Theorem 4.7 for the case $j=4$ (which
appears in \S 7.3).  Indeed, Theorem 4.7 is false as stated.

The two errors are related.  The correction consists of adding the following
hypothesis to the definition (3.1) of a weakly
hyperbolic neighborhood $W$, \begin{itemize}
\item[(d)] $c_1^{-1}(F') \cap c_2^{-1}(F') \cap
W = F$ \end{itemize}
where $F' = c_1(F)= c_2(F)$ as in \S 3.1.  This condition holds
automatically if the correspondence $C$ is the graph of a function.

The addition of condition (d) makes Lemma 5.10 obvious.  Lemma 5.10
is used in \S 5.13 (the proof of Proposition 5.7) which in turn
is needed for the proof of Theorem 4.7.  
The error in \S 7.3 occurs in \S 7.3 Step 2 (the construction of a
morphism $\Phi_4: c_2^*\mathbf{A_4}^{\bullet} \to c_1^!
\mathbf{A_4}^{\bullet}$).  This step should be replaced by the following
construction, which also uses condition (d) above:

Consider the following diagram of spaces and mappings in which
$\diamondsuit$ denotes a
Cartesian square.  Recall that $W\subset C$ is a
neighborhood of the fixed point set $F$ and that $F' = c_1(F) = c_2(F).$
Write
$\widetilde{R}_{Lk} = c_k^{-1}t^{-1}(R_L)$ for $k=1,2$ and define
\[\widetilde{F}_k = F' \times_{R_L} (\widetilde{R}_{L1} \cap \widetilde{R}_{L2})\] 
to be the fiber product in the lower left hand square.  (The mapping $a:t^{-1}(R_L) \to 
t^{-1}(R_L)$ in the lower middle square is the identity.)  The top row is obtained from the
middle row by intersecting with $W$.

\begin{equation*}\begin{CD}
F &@>{h}>>& R_{L1} \cap R_{L2} &@>{a_k}>>& R_{Lk} &@>{j_k}>>& W\\
@V{i_k}VV && @VV{i}V \diamondsuit&& @VV{i_k}V \diamondsuit&&
@VV{i}V \\
\widetilde{F_k} &&\heartsuit& & \widetilde{R}_{L1} \cap
\widetilde{R}_{L2}
&@>{\tilde a_k}>>& \widetilde{R}_{Lk} &@>{\tilde j_k}>>&  C \\
@V{c_k}VV  && @VV{c_k}V && @VV{c_k}V \diamondsuit && @VV{c_k}V
\\
F' &@>>{h_L}>& t^{-1}(R_L) &@>>{a}>& t^{-1}(R_L) &@>>{j_L}>& X
\end{CD}\end{equation*}

We claim that the left-hand rectangle, denoted $\heartsuit$, is Cartesian when $k=1,$
that is, $c_1^{-1}(F') \cap c_1^{-1}t^{-1}(R_L)
\cap c_2^{-1}t^{-1}(R_L) \cap W = F.$   This follows from the hyperbolic hypothesis
together with the (newly added) hypothesis (d) above. 
  
We now  construct a morphism $\phi_4: c_2^*(\mathbf{A}^{\bullet}_4) \to
c_1^!(\mathbf{A}^{\bullet}_4)$ where
$\mathbf{A}^{\bullet}_4 =
j_{L*}h_{L*}h_L^*j_L^{!}\mathbf A^{\bullet}:$
\newcommand{\A}{\mathbf{A}^{\bullet}}

\begin{equation*}
\begin{CD}
c_2^*j_{L*}a_*h_{L*}h_L^*a^*j_L^!\A @>{(2.2)}>>
i_*(i^*c_2^*)j_{L*}a_*h_{L*}h_L^*a^*j_L^!\A @>{(2.5a)}>> \\
i_*j_{2*}(i_2^*c_2^*)a_*h_{L*}h_L^*a^*j_L^!\A @>{(2.2)}>>  
i_*j_{2*}i_2^*\tilde a_{2*}\tilde a_2^*c_2^*a_*h_{L*}h_L^*a^*j_L^! \A
@>{=}>> \\
i_*j_{2*}i_2^*\tilde a_{2*}c_2^*h_{L*}h_L^*a^*j_L^! \A @>{(2.5a)}>>
i_*j_{2*}a_{2*}(i^*c_2^*)h_{L*}h_L^*a^*j_L^!\A @>{(2.5a)}>> \\
i_*j_{2*}a_{2*}h_*i_2^*c_2^*h_L^*a^*j_L^!\A @>{=}>> 
i_*j_{2*}a_{2*}h_*h^*a_2^*i_2^*c_2^*j_L^!\A @>{(2.6a)}>> \\
i_*[j_{2*}a_{2*}h_*h^*a_2^*j_2^!i^*c_2^* \A]
\end{CD} \end{equation*}
Note that $a_2^* = a_2^!$ and $i^* = i^!$ since both $a_2$ and $i$ are
open mappings.  Also, $j_2a_2 = j_1a_1.$  Since $a_1$ is a closed
embedding, $a_{1*} = a_{1!}.$  The above sheaf
$[j_{2*}a_{2*}h_*h^*a_2^*j_2^!i^*c_2^*
\A]$ (in square brackets) is supported on $F$ which is compact and is
contained in $W$, so $i_*[\cdot] = i_![\cdot].$  We also have a morphism
$c_2^*\A \to c_1^!\A.$ Putting these together gives a morphism to
$ i_![j_{1*}a_{1!}h_*h^*a_2^!j_2^!i^!c_1^! \A]$ and hence to
\begin{equation*}\begin{CD}
i_![j_{1*}a_{1!}h_*h^*a_1^!j_1^!i^!c_1^! \A] @>{=}>>
i_!j_{1*}a_{1!}h_*h^*(i^!c_1^!)a^!j_L^!\A @>{(2.6a)}>> \\
i_!j_{1*}a_{1!}h_*(i_1^!c_1^!)h_L^*a^!j_L^!\A @>{(2.4a)}>>
i_!j_{1*}a_{1!}(i^!c_1^!)h_{L*}h_L^*a^!j_L^!\A@>{(2.5b)}>> \\
i_!j_{1*}(i^!c_1^!)a_!h_{L*}h_L^*a^!j_L^!\A @>{(2.4a)}>>
i_!i^!c_1^!j_{L*}a_!h_{L*}h_L^*a^!j_L^!\A @>{(2.2)}>> \\
c_1^!j_{L*}a_!h_{L*}h_L^*a^!j_L^!\A @>{=}>>
c_1^!j_{L*}h_{L*}h_L^*j_L^!\A
\end{CD}\end{equation*}
as desired.  \qed

\subsubsection*{Remarks} Let $\overline{R}(x_0,y_0) = [0,x_0] \times [0,y_0]$ 
denote the closure of the rectangle defined in \S 4.1. 
The new condition (d) implies the following slightly weaker
condition:\begin{itemize}
\item[(d${}^{\prime}$)]
 For some (and hence for any sufficiently small)
$x_0, y_0 > 0$  the set
\[S=  W \cap c_1^{-1}t^{-1}(\overline{R}(x_0,y_0)) \cap
c_2^{-1} t^{-1}(\overline{R}(x_0,y_0))
\subset C\]
is closed (and hence compact),  \end{itemize}
which is equivalent to the statement that the closure of $S$ is contained
in $W.$   The (erroneous) proof of
Lemma 5.10, which appears in \S 5.11, actually shows that (d) implies
(d${}^{\prime}$).   In fact, this
weaker condition (d${}^{\prime}$) suffices for the case $j=1$ of Theorem
4.7, although the stronger condition (d) is needed for the other cases.

\end{document}